\documentclass[12pt]{amsart}\usepackage[numbers,sort&compress]{natbib}
\usepackage{amsmath,amsfonts,amssymb,array,delarray,xspace,a4wide}

\newtheorem{thm}{Theorem}[section]

\newtheorem{prop}[thm]{Proposition}
\newtheorem{cor}[thm]{Corollary}
\newtheorem{conj}[thm]{Conjecture}
\theoremstyle{definition}
\newtheorem{defn}[thm]{Definition}
\newtheorem{rem}[thm]{Remark}
\newtheorem{exmp}[thm]{Example}
\numberwithin{equation}{section}

\newcommand{\mbf}[1]{\ensuremath{\mathbf{#1}}}
\newcommand{\x}{\mbf{x}}\newcommand{\y}{\mbf{y}}
\newcommand{\A}{\mbf{A}}
\newcommand{\one}{\mbf{1}}\newcommand{\zer}{\mbf{0}}
\newcommand{\bbb}[1]{\ensuremath{\mathbb{#1}}}
\newcommand{\N}{\bbb{N}}\newcommand{\Z}{\bbb{Z}}
\newcommand{\R}{\bbb{R}}
\newcommand{\F}{\bbb{F}}
\newcommand{\Zo}{\ensuremath{(\Z^n \perp \one)}}
\newcommand{\Zt}{\ensuremath{(\Z_2^n \perp \one)}}
\newcommand{\smc}[1]{\textrm{\textsc{#1}}}
\newcommand{\sym}{\smc{Sym}}\newcommand{\alt}{\smc{Alt}}
\newcommand{\psl}{\smc{SL}}
\newcommand{\fib}{\smc{Fib}}\newcommand{\luc}{\smc{Luc}}
\newcommand{\lucgp}{\smc{LG}}\newcommand{\fibgp}{\smc{FG}}
\newcommand{\cox}{\smc{Cox}}
\newcommand{\sfw}[1]{\textsf{#1}\xspace}
\newcommand{\sfws}[1]{\textsf{#1}s\xspace}
\newcommand{\ca}{\sfw{CA}}
\newcommand{\sds}{\sfw{SDS}}\newcommand{\sdss}{\sfws{SDS}}
\newcommand{\aca}{\sfw{ACA}}\newcommand{\acas}{\sfws{ACA}}

\newcommand{\Wolf}{\sfw{Wolf}}
\newcommand{\Circ}{\sfw{Circ}}
\newcommand{\per}{\mathsf{Per}}\newcommand{\Fix}{\mathsf{Fix}}
\newcommand{\no}{\mathsf{N}}
\newcommand{\dg}{\mathsf{DG}}\newcommand{\rdg}{\mathsf{RDG}}

\newcommand{\Fy}{\mathfrak{F}_Y}\newcommand{\WOLF}{\mathfrak{Wolf}}
\newcommand{\card}[1]{|#1|}

\newcommand{\floor}[1] {\left\lfloor #1 \right\rfloor}

\begin{document}
%%%%%%%%%%%%%%%%

\title[Dynamics of \acas]{Dynamics groups of asynchronous cellular
  automata}

\author{M.~Macauley}
\author{J.~McCammond}
\author{H.~S.~Mortveit}

\date{May 4, 2010}

\begin{abstract}
  We say that a finite asynchronous cellular automaton (or more
  generally, any sequential dynamical system) is $\pi$-independent if
  its set of periodic points are independent of the order that the
  local functions are applied. In this case, the local functions
  permute the periodic points, and these permutations generate the
  \emph{dynamics group}. We have previously shown that exactly $104$ of
  the possible $2^{2^3}=256$ cellular automaton rules are
  $\pi$-independent. In this article, we classify the periodic states
  of these systems and describe their dynamics groups, which are
  quotients of Coxeter groups. The dynamics groups provide information
  about permissible dynamics as a function of update sequence and, as
  such, connect discrete dynamical systems, group theory, and
  algebraic combinatorics in a new and interesting way. We conclude
  with a discussion of numerous open problems and directions for
  future research.
\end{abstract}



\maketitle

%\tableofcontents

A cellular automaton, or \ca, is a classical discrete dynamical system
defined over a regular grid of cells, such as the lattice $\Z^d$, or
$\Z_n^d$ in the finite case. Every cell takes on one of a finite
number of states and has an update rule that only depends on its state
and the states of its neighbors. Traditionally, at every discrete time
step, the update rules are simultaneously applied. In this article we
study finite cellular automata whose update rules are applied
asynchronously.

More recent work has investigated \emph{sequential dynamical systems}
(\sdss) defined over arbitrary finite graphs where the update rules
are applied asynchronously. (See \cite{Mortveit:07} for a detailed
bibliography.) The asynchronous systems in this paper can be viewed as
either a special type of an \sds or as a modified version of a
classical elementary \ca.
The base graph is the circle graph $\Circ_n$, there are only two
vertex states ($0$ and $1$), and the local update rules are all the
same. There are $2^{2^3}=256$ possible update rules, and in an earlier
article (\cite{Macauley:08a}) we proved that exactly $104$ of these
give rise to an \aca whose periodic states are independent of the
update order, a property independent of the size of the underlying
graph. In this article we describe the periodic states and the
dynamics group for each of these $104$ rules.

If an \aca (or more generally, a sequential dynamical system) is
$\pi$-independent we may construct its \emph{dynamics group}, which is
a permutation group on the set of periodic points. The elements of
this group capture essential information about possible periodic orbit
structures as a function of the update sequence. As such, the group
can be used to characterize the possible long-term dynamics (e.g.,
which periodic orbit configurations can be realized) and also forms an
entry point to the study of permissible long-term dynamics for update
sequence stochastic \acas. Algebraic properties of finite cellular
automata were studied in~\cite{Martin:84}. The introduction of the
dynamics group constitutes a new and promising connection between
algebra and the theory of discrete dynamical systems.

The article is structured as follows. Sections~\ref{sec:sds},
\ref{sec:aca} and~\ref{sec:groups} contain background definitions,
notations, and results. Section~\ref{sec:trivial} studies \acas whose
dynamics groups are trivial, Section~\ref{sec:invertible} focuses on
rules that are invertible, and Section~\ref{sec:exceptional}
investigates the remaining cases. The results are summarized in
Tables~\ref{tbl:trivial},~\ref{tbl:invertible},
and~\ref{tbl:exceptional}. We conclude with a discussion of future
research problems and potential applications to those working on
cellular automata.

%%%%%%%%%%%%%%%%%%%%%%%%%%%%%%%%%%%%%%%%%%%%%%%%%%%%%
\section{Sequential Dynamical Systems}\label{sec:sds}

We will present the main concepts in this article, such as the
dynamics group, in the more general setting of sequential dynamics
systems. A \emph{sequential dynamical system}, or \sds, is a discrete
dynamical system with three components: an undirected graph $Y$, a
list of update rules $\Fy$, and an update order $\omega$. Seeing the
setup in its full generality is not only useful, but it is consistent
with the notation in the prequel to this paper~\cite{Macauley:08a}
that contains the classification of the $104$ $\pi$-independent \acas,
and the original paper on dynamics groups~\cite{Hansson:05b}.

\begin{defn}[Graph conventions]
  Let $Y$ be a simple undirected graph with $n$ vertices labeled from
  $1$ to $n$, and recall that the \emph{neighbors} of a vertex are
  those vertices connected to it by an edge.  If $\F$ is a finite
  field and every vertex is assigned a value from $\F$, then a global
  state of the system is described by an $n$-tuple $\y$ whose $i^{\rm
    th}$ coordinate indicates the current state of the vertex $i$.
  The set of all possible states is the vector space $\F^n$.
\end{defn}

\begin{defn}[Local functions]
  A function $F\colon\F^n\to \F^n$ is called \emph{$Y$-local at $i$}
  if for each $\y\in \F^n$ (1) $F(\y)$ only alters the $i^{\rm th}$
  coordinate of $\y$ and (2) the new value of the $i^{\rm th}$
  coordinate only depends on the coordinates of $\y$ corresponding to
  $i$ and its neighbors in $Y$.  Other names for such a function are a
  \emph{local function} or an \emph{update rule}.  We use $\Fy$ to
  denote a list with one local function for each vertex of $Y$.  More
  precisely, $\Fy = (F_1, F_2, \ldots, F_n)$ where $F_i$ is a function
  that is $Y$-local at $i$.
\end{defn}

\begin{defn}[Restricted local functions]
  If $i$ is a vertex with $k$ neighbors in $Y$, then corresponding to
  each function $F$ that is $Y$-local at $i$, we define a function
  $f\colon\F^{k+1}\to \F$ where the domain is restricted to the
  coordinates corresponding to $i$ and its neighbors, and the output
  is the new value $F$ would assign to the $i^{\rm th}$ coordinate
  under these conditions.  The functions $F$ and $f$ contain the same
  information packaged differently and each determines the other.
  Both have their uses.  Functions such as $F$ can be readily
  composed, but functions such as $f$ are easier to explicitly
  describe.
\end{defn}

\begin{defn}[Update orders]
  An \emph{update order} $\omega$ is a finite sequence of numbers
  chosen from the set $\{1,\ldots, n\}$.  If every number $1 \leq i
  \leq n$ occurs at least once, we say it is \emph{fair} and if every
  number occurs exactly once, then it is \emph{simple}.  We use the
  notation $\omega = (\omega_1, \omega_2, \ldots, \omega_m)$ with
  $m=|\omega|$.  Of course, $m \geq n$ when $\omega$ is fair and $m=n$
  when $\omega$ is simple.  Let $W_Y$ denote the collection of all
  update orders and let $S_Y$ denote the subset of simple update
  orders.  The subscript $Y$ indicates that we are thinking of the
  numbers in these sequences as vertices in the graph $Y$.
\end{defn}

\begin{defn}[Sequential dynamical systems]
  A \emph{sequential dynamical system}, or \sds, is a triple
  $(Y,\Fy,\omega)$ consisting of an undirected graph $Y$, a list of
  local functions $\Fy$, and a fair update order $\omega\in W_Y$.  If
  $\omega$ is the sequence $(\omega_1,\omega_2,\ldots,\omega_m)$, then
  we construct the \emph{\sds map} $[\Fy,\omega]\colon \F^n \to \F^n$
  as the composition $[\Fy,\omega]:=F_{\omega_m} \circ \cdots \circ
  F_{\omega_1}$.
\end{defn}

The main goal is to understand the dynamics of the \sds map, i.e., its
behavior under iteration.

\begin{defn}[$\pi$-independence]
  Let $\per[\Fy,\omega]$ denote the set of states periodic under
  iterations of $[\Fy,\omega]$. A list of $Y$-local functions $\Fy$
  is called \emph{$\omega$-independent} if $\per[\Fy,\omega] =
  \per[\Fy,\omega']$ for all fair update orders $\omega, \omega' \in
  W_Y$ and \emph{$\pi$-independent} if $\per[\Fy,\pi] =
  \per[\Fy,\pi']$ for all simple update orders $\pi, \pi' \in
  S_Y$. When $\Fy$ is $\pi$-independent, we write $\per(\Fy)$
  instead of $\per[\Fy,\pi]$.
\end{defn}

In the case of $\pi$-independence we may -- by abuse of notation --
let the list of $Y$-local functions $\Fy$ stand for the entire
\sds. When $\Fy$ is $\pi$-independent, its local functions permute the
elements of $\per(\Fy)$.

\begin{prop}[Permuting periodic states]\label{prop:bij}
  If $\Fy$ is $\pi$-independent and $P=\per(\Fy)$, then for each
  $i$, $F_i(P) = P$.  In particular, the restriction of $F_i$ to $P$
  is a permutation.
\end{prop}

\begin{proof}
  Let $\omega=(\pi_1,\pi_2,\dots,\pi_n)$ be a simple update
  order with $\pi_1 =i$ and let $\sigma$ be the modified update
  order with $\pi_1$ moved from the first to last:
  $\sigma=(\pi_2,\pi_3,\dots,\pi_n,\pi_1)$.  Since
  $F_i\circ[\Fy,\pi]^k = [\Fy,\sigma]^k\circ F_i$ for all $k$, and
  by hypothesis $[\Fy,\pi]^k(\F^n) = [\Fy,\sigma]^k(\F^n)=P$ for
  all sufficiently large $k$, we find that $F_i(P)\subset P$.  More
  explicitly, for large enough $k$,
  \[
  F_i(P)=F_i \circ [\Fy,\pi]^k(\F^n)= [\Fy,\sigma]^k \circ
  F_i(\F^n) \subset P.
  \]
  Moreover, since $[\Fy,\pi](P)=P$, the restriction of $F_i$ to $P$
  is injective, and $F_i(P)=P$.
\end{proof}

When $\Fy$ is $\pi$-independent, we write $F_i^*$ and $[\Fy,\pi]^*$ to
denote the restrictions of these maps to $\per(\Fy)$.  By
Proposition~\ref{prop:bij} all such maps are permutations.  Note that
$\pi$-independence focuses on the periodic states as a set rather than
how these states are permuted.  In particular, when a
$\pi$-independent $\Fy$ is paired with two different update orders
$\pi$ and $\sigma$, the permutations $[\Fy,\pi]^*$ and
$[\Fy,\sigma]^*$ are often distinct.  These various permutations can
be used to construct a group encoding all of the possible
dynamics~\cite{Mortveit:07}.

\begin{defn}[Dynamics group]
  Let $\Fy$ be $\pi$-independent.  For any collection of update
  orders $U\subseteq W_Y$, the {\em dynamics group of $\Fy$ with
    respect to $U$} is
%%the group generated by the corresponding
%%  permutations of $\per(\Fy)$.
\begin{equation*}
  \dg(\Fy,U) = \langle [\Fy,\omega]^* \mid \omega\in U \rangle \;.
\end{equation*}
%%
%% Our notation for this group is $\dg(\Fy,U)$.
  It should be clear that when $U$ and $V$ are sets of update orders
  and $U$ is contained in the closure of $V$ under concatenation, then
  $\dg(\Fy,U) \subset \dg(\Fy,V)$.  The {\em dynamics group} of $\Fy$,
  $\dg(\Fy)=\dg(\Fy,W_Y)$, and the {\em restricted dynamics group} of
  $\Fy$, $\rdg(\Fy)=\dg(\Fy,S_Y)$, are special cases of particular
  interest.  Note that $\dg(\Fy)$ contains and is generated by the
  bijections $F_i^*$.
\end{defn}

%%%%%%%%%%%%%%%%%%%%%%%%%%%%%%%%%%%%%%%%%%%%%%%%%%%%%%%
\section{Asynchronous Cellular Automata}\label{sec:aca}

The \sdss we focus on are defined over circular graphs, they have only
two possible vertex states, and all the local functions are
identically defined.  These \emph{asynchronous cellular automata}, or
\acas, are an asynchronous version of the classical finite elementary
cellular automata.  Even in such a restrictive situation there are
many interesting dynamical behaviors.

\begin{defn}[Circular graphs and vertex states]
  Let $Y=\Circ_n$ denote the \emph{circular graph} with vertex set
  $\{1,\dots,n\}$ (viewed as residue classes mod $n$) and edges
  connecting $i$ and $i+1$ mod $n$.  To avoid trivialities, we always
  assume $n>3$.  Each vertex has two possible states that we identify
  with $\F_2 = \{0,1\}$, the field of size $2$.
\end{defn}

\begin{defn}[Wolfram rules]
  Let $F_i\colon \F_2^n\to \F_2^n$ be a function $\Circ_n$-local at
  $i$ and let $f_i\colon \F_2^3\to \F_2$ be its restricted form.
  Because the neighbors of $i$ are $i-1$ and $i+1$, it is conventional
  to list these coordinates in ascending order in the domain of $f_i$,
  keeping in mind that all subscripts are viewed mod $n$.  The
  function $F_i$, henceforth referred to as a \emph{Wolfram rule},
  updates the value of $y_i$ based on the value of the triple
  $(y_{i-1},y_i,y_{i+1})$ and it is completely determined by how the
  $i^{\rm th}$ coordinate is updated in these eight possible
  situations.  In other words, $F_i$ is completely described by the
  following table.
  \[
    \begin{array}{c||c|c|c|c|c|c|c|c}
      y_{i-1}y_iy_{i+1} & 111 & 110 & 101 & 100 & 011 & 010 & 001 & 000
      \\ \hline f_i(y_{i-1},y_i,y_{i+1}) & a_7 & a_6 & a_5 & a_4 & a_3
      & a_2 & a_1 & a_0
    \end{array}
  \]
  More concisely, the $2^8=256$ possible Wolfram rules can be indexed
  by an $8$-digit binary number $a_7a_6a_5a_4a_3a_2a_1a_0$, or by its
  decimal equivalent $k = \sum a_i 2^i$.  There is thus a
  \emph{Wolfram rule $k$} for each integer $0\leq k \leq 255$.
\end{defn}

\begin{defn}[Asynchronous cellular automata]
  We write $\Wolf^{(k)}_i$ to denote the update rule $F_i\colon
  \F^n\to \F^n$ corresponding to $k$ and $\WOLF^{(k)}_n$ for the list
  $(\Wolf^{(k)}_1, \Wolf^{(k)}_2, \ldots, \Wolf^{(k)}_n)$ of update
  rules.  For each fair update order $\omega$ the \sds
  $(\Circ_n,\WOLF_n^{(k)},\omega)$ is called an {\em asynchronous
  cellular automaton}, or \aca.  If $\WOLF_n^{(k)}$ is
  $\pi$-independent ($\omega$-independent) for all $n>3$, we say
  Wolfram rule~$k$ is $\pi$-independent ($\omega$-independent).
\end{defn}

When $\WOLF_n^{(k)}$ is $\pi$-independent, let $P_{n,k} =
\per(\WOLF_n^{(k)})$ denote its periodic states and let $G_{n,k} =
\dg(\WOLF_n^{(k)})$ denote its dynamics group. We usually suppress the
dependence on $n$ and write simply $P_k$ and $G_k$.  In this notation,
our goal is to describe the set $P_k$ and the group $G_k$ for each
$\pi$-independent Wolfram rule. In~\cite{Macauley:08a} we proved the
following result.

\begin{thm}\label{thm:104}
  Exactly $104$ Wolfram rules are $\pi$-independent.  More
  precisely, $\WOLF_n^{(k)}$ is $\pi$-independent for all $n>3$ iff
  $k \in \{0$, $1$, $4$, $5$, $8$, $9$, $12$, $13$, $28$, $29$, $32$,
  $40$, $51$, $54$, $57$, $60$, $64$, $65$, $68$, $69$, $70$, $71$,
  $72$, $73$, $76$, $77$, $78$, $79$, $92$, $93$, $94$, $95$, $96$,
  $99$, $102$, $105$, $108$, $109$, $110$, $111$, $124$, $125$, $126$,
  $127$, $128$, $129$, $132$, $133$, $136$, $137$, $140$, $141$,
  $147$, $150$, $152$, $153$, $156$, $157$, $160$, $164$, $168$,
  $172$, $184$, $188$, $192$, $193$, $194$, $195$, $196$, $197$,
  $198$, $199$, $200$, $201$, $202$, $204$, $205$, $206$, $207$,
  $216$, $218$, $220$, $221$, $222$, $223$, $224$, $226$, $228$,
  $230$, $232$, $234$, $235$, $236$, $237$, $238$, $239$, $248$,
  $249$, $250$, $251$, $252$, $253$, $254$, $255\}$.
\end{thm}

In~\cite{Macauley:08a} we also defined the inversion, reflection, and
inversion-reflection of an \aca. Loosely speaking, \emph{inversion}
systematically swaps the roles of $0$ and $1$, \emph{reflection}
systematically switches left and right, and
\emph{inversion-reflection} does both at once. The inversion,
reflection, or inversion-reflection of a $\pi$-independent Wolfram
rule is still $\pi$-independent, it has a corresponding set of
periodic states and an isomorphic dynamics group.  This should not be
surprising since all we have done is relabel the underlying states on
which the local functions act.  Rules related in this manner are
%%said to be
{\em dynamically equivalent}.  When the $256$ Wolfram rules are
partitioned into classes of rules related by reflection, inversion or
both, there are $88$ equivalence classes. The $104$ rules listed in
Theorem~\ref{thm:104} belong to $41$ such classes and thus we only
need to describe $P_k$ and $G_k$ for $41$ representative values of
$k$.

\begin{cor}[$41$ representative rules]
  Every $\pi$-independent Wolfram rule is dynamically equivalent to
  $\WOLF_n^{(k)}$ for some $k \in \{0$, $1$, $4$, $5$, $8$, $9$, $12$,
  $13$, $28$, $29$, $32$, $40$, $51$, $54$, $57$, $60$, $72$, $73$,
  $76$, $77$, $105$, $128$, $129$, $132$, $133$, $136$, $137$,
  $140$, $141$, $150$, $152$, $156$, $160$, $164$, $168$, $172$,
  $184$, $200$, $201$, $204$, $232\}$.
\end{cor}

Finally, because the binary notation is cumbersome and the decimal
notation is opaque, we introduce (as in \cite{Macauley:08a}) a concise
symbolic tag for each Wolfram rule.

\begin{defn}[Tags]
  The four functions from $\F_2$ to $\F_2$ can be described by their
  behavior: the value never changes, the value always changes, both
  elements go to $0$, or both elements go to $1$.  We refer to these
  functions by the evocative symbols {\tt -}, {\tt x}, {\tt 0}, and
  {\tt 1}, respectively.  The restricted local form of Wolfram rule
  $k$ is completely determined by the four functions from $\F_2$ to
  $\F_2$ that result when the values of $y_{i-1}$ and $y_{i+1}$ are
  held constant.  Let $t_0$, $t_1$, $t_2$ and $t_3$ be the symbols for
  these functions in the four cases $y_{i-1}=y_{i+1}=0$, $y_{i-1}=0$
  and $y_{i+1}=1$, $y_{i-1}=1$ and $y_{i+1}=0$, and
  $y_{i-1}=y_{i+1}=1$, respectively.  The \emph{tag} of $k$ is the
  string $t_3t_2t_1t_0$.  Note that $t_0$ depends on the values of
  $a_0$ and $a_2$, $t_1$ depends on $a_1$ and $a_3$, $t_2$ depends on
  $a_4$ and $a_6$, and $t_3$ depends on $a_5$ and $a_7$.  The
  numbering and the order of the $t_i$'s has been chosen to match the
  traditional binary representation as closely as possible, easing the
  transition between the two.  As an illustration, the reader can
  verify that Wolfram rule $29$ has binary notation $00011101$ and tag
  {\tt 0x-1}.
\end{defn}

Patterns among the Wolfram rules are easier to discern when using
tags.

\begin{rem}[Tags and dynamic equivalence]
  On the level of tags, reflections switch the order of $t_1$ and
  $t_2$.  For example, the reflection of rule {\tt0-1x} is {\tt 01-x}.
  To describe the effect that inversion has on tags, we define a map
  $\iota \colon \{{\tt 1,0,-,x}\} \to \{{\tt 1,0,-,x}\}$ that fixes
  {\tt -} and {\tt x} while switching {\tt 0} and {\tt 1}.  When
  Wolfram rule $k$ has tag $t_3t_2t_1t_0$, its inversion has tag
  $\iota(t_0) \iota(t_1) \iota(t_2) \iota(t_3)$.  For example, the
  inversion of rule {\tt 0-1x} is {\tt x0-1}.  See \cite{Macauley:08a}
  for a more detailed explanation.
\end{rem}

\begin{rem}[Tags and other \sdss]
  If $(Y,\Fy,\omega)$ is an \sds with only two vertex states, then
  each update rule $F_i$ can be described by a set of symbols similar
  to the tag used to describe Wolfram rules.  More specifically, if
  vertex $i$ has exactly $k$ neighbors, then the restricted local form
  of $F_i$ is a function $f_i:\F_2^{k+1}\to \F_2$ and this function is
  determined by the $2^k$ functions from $\F_2$ to $\F_2$ that result
  when the values of the neighbors of $i$ are held constant.  In
  particular, the behavior of $F_i$ is determined by the corresponding
  $2^k$ symbols, selected from $\{{\tt 1,0,-,x}\}$.
\end{rem}

%%%%%%%%%%%%%%%%%%%%%%%%%%%%%%%%%%
\section{Groups}\label{sec:groups}
%%%%%%%%%%%%%%%%%%%%%%%%%%%%%%%%%%

In preparation for our investigation of dynamics groups we recall a
few basic facts about group actions and Coxeter groups.

\begin{defn}[Faithful actions]
  Let $G$ be a group acting on a set $X$.  The action is called
  \emph{faithful} if nontrivial elements act nontrivially.  When this
  is the case, we can view $G$ as a subgroup of $\sym_X$.  Note that
  groups generated by permutations act faithfully on their underlying
  sets almost by definition.  In particular, for any $\pi$-independent
  \sds $\Fy$, the action of its dynamics group $\dg(\Fy)$ on its
  periodic states $\per(\Fy)$ is faithful.
\end{defn}

One of the key features of a group acting on a set is its orbit
structure.

\begin{defn}[Orbits]
  Let $G$ be a group acting on a finite set $X$.  The \emph{orbit} of
  a point $x\in X$ is the subset of points to which it can be sent by
  an element of $G$.  Thus the orbit of $x$ is $Gx = \{ g \cdot x \mid
  g\in G\}$.  Because two orbits are either identical or disjoint, the
  collection of all orbits $\{ Gx \mid x\in X \}$ partitions $X$ into
  equivalence classes $X = X_1 \sqcup X_2 \sqcup \ldots \sqcup
  X_\ell$.  When there is only one orbit we say the action is
  \emph{transitive}.
\end{defn}

As we investigate the faithful action of $G_{n,k}$ on $P_{n,k}$, we
use $o_{n,k}$ to denote the number of orbits under this action.  And,
as with $G_k$ and $P_k$, we usually suppress the dependence on $n$ and
write $o_k$ instead.  The orbit structure is of interest because
distinct orbits provide information about the structure of the group.

\begin{defn}[Subdirect products]
  If $G$ acts faithfully on a finite set $X$ with orbits $X_1$, $X_2$,
  \ldots, $X_\ell$ then $G$ can be viewed as a subgroup of $\sym_{X_1}
  \times \sym_{X_2} \times \cdots \times \sym_{X_{\ell}}$.  To see
  this note that when an element of $G$ is written in disjoint cycle
  notation, each cycle must permute elements within a single orbit.
  Thus every $g \in G$ can be viewed as a $\ell$-tuple $g =
  (g_1,g_2,\ldots, g_\ell)$ where $g_i \in \sym(X_i)$.  Even better,
  we can replace each $\sym_{X_i}$ with $G_i$, the image of $G$ under
  the projection to the $i^{th}$ factor.  As a result $G$ embeds in a
  direct product $G_1 \times G_2 \times \cdots \times G_\ell$, where
  the projection to each factor is onto.  We call $G$ a
  \emph{subdirect product} of $G_1$, $G_2$,\ldots, and $G_\ell$, and
  the $G_i$s are called the \emph{orbit groups} of $G$.
\end{defn}

For an illustration of these concepts, consider a group generated by a
single permutation.

\begin{exmp}[Cyclic groups]
  The permutation $(1,2)(3,4,5)(6,7,8,9)$ generates a cyclic group $G$
  of order $12$.  It naively belongs to $\sym_9$, a group of size
  $362880$, but based on its orbit structure it lives inside the much
  smaller group $\sym_2\times \sym_3 \times \sym_4$ of size $288$.
  And if we cut down each factor to the image under projection, then
  $G$ embeds in $\Z_2 \times \Z_3 \times \Z_4$, a noncyclic group of
  order $24$ that contains $G$ as an index~$2$ subgroup.
\end{exmp}

One caution is that $G$ can be a subdirect product of groups without
splitting as a direct product.  Consider the cyclic group $G$
generated by $(1,2,3)(4,5,6)$.  The procedure described above embeds
$G$ as a subdirect product of $\Z_3$ and $\Z_3$, but $G$ itself is
simple.
The second set of results we need to recall are about Coxeter groups.
A Coxeter group is a group generated by involutions with a
presentation of a particularly simple form.

\begin{defn}[Coxeter groups]
  Let $W$ be a group generated by a finite subset $S =
  \{s_1,\ldots,s_n\}$ and let $m_{ij} \in \N \cup \{\infty\}$ denote
  the order of the product $s_is_j$.  The pair $(W,S)$ is called a
  \emph{Coxeter system} if the elements of $S$ are involutions
  (i.e., $m_{ii} =1$) and the presentation $\langle
  s_1,\dots,s_n\,|\,(s_is_j)^{m_{ij}}=1\rangle$ is a presentation of
  $W$.  Note that when $m_{ij} = \infty$, no relation is included for
  this pair of indices and that trivially $m_{ij}=m_{ji}\geq 2$ for
  all $i\neq j$.
\end{defn}

The letters $W$ and $S$ are those traditionally used for a
\emph{Coxeter group} and its \emph{Coxeter generators}.  Coxeter
groups have an incredibly rich theory and close connections with many
areas of mathematics.  By comparison, the results we need are fairly
modest.  Before listing these results, we first establish the
relevance of Coxeter groups to our investigation of dynamics groups.

\begin{prop}[Coxeter quotients]\label{prop:quotient}
  Every group generated by a finite set of involutions can be viewed
  as a quotient of a Coxeter group in a natural way.
\end{prop}

\begin{proof}
  Let $G$ be a group and let $S = \{s_1,\ldots,s_n\}$ be a subset of
  involutions that generate $G$.  If we define $m_{ij}$ as the order
  of $s_is_j$ in $G$ and we define $W$ as the group with presentation
  $\langle s_1,\dots,s_n\,|\,(s_is_j)^{m_{ij}}=1\rangle$ then there is
  a natural surjective homomorphism from $W$ to $G$ sending $s_i$ to
  $s_i$.
\end{proof}

\begin{thm}[Dynamics groups as Coxeter quotients]
  If $\Fy$ is a $\pi$-independent \sds with only two possible
  vertex states, then $F_i^*$, the restriction of a local function to
  the periodic states $\per(\Fy)$, is either trivial or an involution.
  As a consequence, the dynamics group $\dg(\Fy)$ is either trivial or
  a quotient of a Coxeter group.
\end{thm}

\begin{proof}
  The key observation is that because $F_i^*$ can only change the
  $i^{th}$ coordinate, the size of the cycles in its cycle structure
  are bounded by the number of possible vertex states.  For the second
  assertion, note that when the dynamics group $\dg(\Fy)$ is
  nontrivial, it is generated by the nontrivial $F_i^*$ and then apply
  Proposition~\ref{prop:quotient}.
\end{proof}

The case when $\dg(\Fy)$ is trivial can be recognized by its fixed
points.

\begin{defn}[Fixed points]
  If $\y\in \per(\Fy)$ is fixed under some simple update order $\pi$,
  then $\y$ must be fixed by each local function $F_i$.  This is
  because a change to the $i^{th}$ coordinate cannot be corrected by
  the other local functions in the composition that produces
  $[\Fy,\pi]$.  As a consequence, $\y$ is fixed under \emph{all}
  update orders $\omega$.  We write $\Fix(\Fy)$ to denote the set of
  periodic states fixed by some simple update order, or equivalently,
  the set of periodic states fixed by all local functions $F_i$ (and
  thus fixed by all update orders).
\end{defn}

\begin{prop}[Trivial groups and fixed points]\label{prop:trivial}
  A $\pi$-independent \sds $\Fy$ has a trivial dynamics group
  $\dg(\Fy)$ iff $\Fix(\Fy)=\per(\Fy)$.
\end{prop}

\begin{proof}
  If $\Fix(\Fy)=\per(\Fy)$ then each $F_i^*$ is trivial and $\dg(\Fy)$
  is trivial.  On the other hand, if $\Fix(\Fy)\neq \per(\Fy)$ then
  there is a periodic state $\y$ and a local function $F_i$ such that
  $F_i(\y) \neq \y$.  For this $i$, $F_i^*$ is nontrivial and thus
  $\dg(\Fy)$ is nontrivial.
\end{proof}

When the dynamics group $\dg(\Fy)$ is nontrivial, the graph $Y$ can be
used to describe the Coxeter group of which it is a quotient.

\begin{defn}[Coxeter diagrams]
  The presentation of a Coxeter group is often summarized in graphical
  form as follows.  Given a Coxeter system $(W,S)$ we construct a
  graph with vertices indexed by $S$ and an edge labeled $m_{ij}$
  connecting vertex $i$ and $j$ whenever $m_{ij} >2$.  The Coxeter
  presentation of $W$ can be easily reconstructed from this
  edge-labeled graph called the \emph{Coxeter diagram} of $W$.
\end{defn}

\begin{defn}[Coxeter diagrams for dynamics groups]\label{def:coxeter-dg}
  If $\Fy$ is a $\pi$-independent \sds with only two possible vertex
  states and a nontrivial dynamics group $\dg(\Fy)$, then the Coxeter
  diagram for the Coxeter group of which $\dg(\Fy)$ is a quotient can
  be obtained from $Y$ in three easy steps.  First remove every vertex
  $i$ (and the edges connected to it) for which $F_i^*$ is trivial.
  Next, remove the edges between $i$ and $j$ when $m_{ij}=2$ (or
  equivalently when $F_i^*$ and $F_j^*$ commute).  And finally, add
  the label $m_{ij}$ to each remaining edge.  We note that because
  $F_i^*$ and $F_j^*$ can only alter coordinates $i$ and $j$, the
  cycles of $F_i^*\circ F_j^*$ have size at most $4$.  Thus each
  $m_{ij}$ divides $12$, the gcd of possible cycle lengths.
\end{defn}

\begin{defn}[Coxeter label]\label{def:coxeter-label}
  The Coxeter diagram for the Coxeter group mapping onto the dynamics
  group of a $\pi$-independent Wolfram rule $k$ is particularly simple
  because of the symmetry of construction of the $\WOLF_n^{(k)}$.
  Writing $F_i$ instead of $\Wolf_i^{(k)}$ we see that one $F_i^*$ is
  nontrivial iff they are all nontrivial and the order of $F_i^* \circ
  F_{i+1}^*$ is a constant, independent of $i$ and $n$.  We call this
  constant $c_k$, the \emph{Coxeter label} of Wolfram rule $k$, and as
  noted above, the value of $c_k$ must divide $12$.  When $G_k$ is
  nontrivial all vertices remain.  If $c_k>2$, all edges remain and
  are labeled $c_k$.
\end{defn}

The result we need from Coxeter theory is an identification of certain
classes of groups.

\begin{rem}[Small Coxeter labels]\label{rem:cox}
  If $\WOLF_n^{(k)}$ is $\pi$-independent, $G_k$ is nontrivial, and
  $c_k=2$, then $G_k$ is a quotient of the Coxeter group $\Z_2^n$
  defined by an edgeless Coxeter diagram.  If $\WOLF_n^{(k)}$ is
  $\pi$-independent, $G_k$ is nontrivial, and $c_k=3$, then $G_k$
  is a quotient of the Coxeter group defined by a circular Coxeter
  diagram with edges labeled $3$.  This group is the affine Coxeter
  group of type $\widetilde A_{n-1}$ but since it is the only Coxeter
  group we consider without a pre-existing common name (such as
  $\sym_n$), we call this group $\cox_n$.  Its structure is
  well-known.  Let $\Zo$ denote the subset of $\Z^n$ perpendicular to
  the vector $\one = (1,1,\ldots,1)$, i.e., the set of vectors whose
  coordinates sum to $0$, and note that these form a subgroup under
  vector addition.  If we let $\sym_n$ act on $\Zo$ by permuting
  coordinates in the natural way, then $\cox_n$ is isomorphic to the
  semidirect product $\Zo \rtimes \sym_n$.  Geometrically, it is the
  group of isometries of the Euclidean space $\R^n$ that preserve the
  sum of the coordinates and send vectors with all integer coordinates
  to other such vectors.

  For the sake of concreteness, we select the following explicit
  isomorphism between $\cox_n$ and $\Zo \rtimes \sym_n$.  To avoid
  confusion we use $\x = (x_1,x_2,\ldots,x_n)$ for an element of $\Zo$
  and reserve $\y = (y_1,y_2,\ldots, y_n)$ for states.  Let
  $s_1,s_2,\ldots,s_n$, denote the Coxeter generators of $\cox_n$.
  For each $i<n$ let $s_i$ be the element that switches $x_i$ and
  $x_{i+1}$, and let $s_n$ send $\x = (x_1,x_2,\ldots, x_{n-1}, x_n)$
  to $(x_n-1,x_2,\ldots,x_{n-1},x_1+1)$.  The natural map to $\sym_n$
  only remembers how the subscripts on the $x_i$s are permuted and
  ignores the translational component.  An element in the kernel of
  the map $\cox_n \to \sym_n$ is called a \emph{pure translation}
  since it merely adds to $\x$ a vector in $\Zo$.  These pure
  translations are generated by the elements $T_i$ in $\cox_n$ defined
  as follows.  The element $T_1= s_1 s_2 \cdots s_{n-1} s_n s_{n-1}
  \cdots s_3 s_2$ and the other $T_i$ are obtained by consistently
  shifting the subscripts.  In terms of its effect on an element
  $\x\in \Zo$, $T_1$ adds the vector $\langle
  -1,0,\ldots,0,1\rangle$. More generally, $T_i$ adds $1$ to $x_{i-1}$
  and subtracts $1$ from $x_i$.  The pure translations are generated
  by these commuting $T_i$ which are nearly independent. The sole
  nontrivial relation they satisfy is that $T_1+T_2+\cdots T_n$ is
  trivial.
\end{rem}

%%%%%%%%%%%%%%%%%%%%%%%%%%%%%%%%%%%%%%%%%%%
\section{Trivial groups}\label{sec:trivial}
%%%%%%%%%%%%%%%%%%%%%%%%%%%%%%%%%%%%%%%%%%%

In this section we discuss $\pi$-independent Wolfram rules with
trivial dynamics group.  Of the $41$ representative Wolfram rules,
$26$ of them fall into this category.  For rule $204$ (with tag {\tt
  ----}), this is immediate since as its tag indicates no local
function ever alters the current state.  For the other $25$ rules the
triviality of $G_k$ is a consequence of
Proposition~\ref{prop:trivial}.  More specifically, several lemmas in
\cite{Macauley:08a} established conditions under which all periodic
points are fixed.  These lemmas are listed below along with the
representative rules they cover.

\begin{itemize}
\item Lemma~5.3 (Rules $0$, $4$, $8$, $12$, $72$, $76$, $128$, $132$, $136$,
$140$ and $200$)
\item Lemma~5.5 (Rules $160$, $164$, $168$, $172$ and $232$)
\item Lemma~5.6 (Rules $5$, $13$, $77$, $133$ and $141$)
\item Lemma~6.1 (Rules $32$ and $40$)
\item Lemma~6.3 (Rules $152$ and $184$)
\end{itemize}

Even though the group $G_k$ is trivial in each case, the set $P_k$
remains to be calculated.  (Many of the proofs in~\cite{Macauley:08a}
are nonconstructive and do not determine the set of periodic states
explicitly.)  Because periodic states and fixed states coincide for
these rules, there is a relatively straightforward procedure for
finding them: simply look at the definition of the rule and remove all
states containing triples $y_{i-1}y_iy_{i+1}$ of consecutive states
that would lead to an alteration.  Clearly the removed states are not
fixed by all local functions and do not belong to $P_k$, and any
states remaining at the end of this procedure are fixed by every local
function and do belong in $P_k$.  This prompts the following
definitions.

\begin{table}
\begin{tabular}{|c|c|rrr|c|c|c|c|}
  \hline
  rule & tag & inv & refl& i.r. & $o_k$ & $P_k$ & $c_k$ & $G_k$ \\  \hline  \hline
   0 & \tt{0000} & 255 &   0 & 255 & $1$ & $\zer$ & $1$ & $1$ \\
   8 & \tt{00-0} & 239 &  64 & 253 & $1$ & $\zer$ & $1$ & $1$ \\
  32 & \tt{x000} & 251 &  32 & 251 & $1$ & $\zer$ & $1$ & $1$ \\
  40 & \tt{x0-0} & 235 &  96 & 249 & $1$ & $\zer$ & $1$ & $1$ \\
 \hline
 128 & \tt{-000} & 254 & 128 & 254 & $2$ & $\zer\cup\one$ & $1$ & $1$ \\
 136 & \tt{-0-0} & 238 & 192 & 252 & $2$ & $\zer\cup\one$ & $1$ & $1$ \\
 160 & \tt{1000} & 250 & 160 & 250 & $2$ & $\zer\cup\one$ & $1$ & $1$ \\
 168 & \tt{10-0} & 234 & 224 & 248 & $2$ & $\zer\cup\one$ & $1$ & $1$ \\
 \hline
 152 & \tt{-x-0} & 230 & 194 & 188 & $2$ & $\zer\cup\one$ & $1$ & $1$ \\
 184 & \tt{1x-0} & 226 & 226 & 184 & $2$ & $\zer\cup\one$ & $1$ & $1$ \\
 \hline
   4 & \tt{000-} & 223 &   4 & 223 & $\card{\no_{A}}$ & $\no_{A}$ & $1$ & $1$ \\
  12 & \tt{00--} & 207 &  68 & 221 & $\card{\no_{A}}$ & $\no_{A}$ & $1$ & $1$ \\
 132 & \tt{-00-} & 222 & 132 & 222 & $\card{\no_{A}}+1$ & $\no_{A}\cup\one$ & $1$ & $1$ \\
 140 & \tt{-0--} & 206 & 196 & 220 & $\card{\no_{A}}+1$ & $\no_{A}\cup\one$ & $1$ & $1$ \\
 \hline
   5 & \tt{0001} &  95 &   5 &  95 & $\card{\no_{AB}}$ & $\no_{AB}$ & $1$ & $1$ \\
  13 & \tt{00-1} &  79 &  69 &  93 & $\card{\no_{AB}}$ & $\no_{AB}$ & $1$ & $1$ \\
 133 & \tt{-001} &  94 & 133 &  94 & $\card{\no_{AB}}+1$ & $\no_{AB}\cup\one$ & $1$ & $1$ \\
 141 & \tt{-0-1} &  78 &  197 & 92 & $\card{\no_{AB}}+1$ & $\no_{AB}\cup\one$ & $1$ & $1$\\
 \hline
 164 & \tt{100-} & 218 & 164 & 218 & $\card{\no_{AE}}+1$ & $\no_{AE}\cup\one$ & $1$ & $1$ \\
 172 & \tt{10--} & 202 & 228 & 216 & $\card{\no_{AE}}+1$ & $\no_{AE}\cup\one$ & $1$ & $1$ \\
 \hline
  77 & \tt{0--1} &  77 &  77 &  77 & $\card{\no_{BC}}$ & $\no_{BC}$ & $1$ & $1$ \\
  76 & \tt{0---} & 205 &  76 & 205 & $\card{\no_{C}}$ & $\no_{C}$ & $1$ & $1$ \\
  72 & \tt{0--0} & 237 &  72 & 237 & $\card{\no_{CD}}$ & $\no_{CD}$ & $1$ & $1$ \\
 200 & \tt{---0} & 236 & 200 & 236 & $\card{\no_{D}}$ & $\no_{D}$ & $1$ & $1$ \\
 232 & \tt{1--0} & 232 & 232 & 232 & $\card{\no_{DE}}$ & $\no_{DE}$ & $1$ & $1$ \\
 \hline
 204 & \tt{----} & 204 & 204 & 204 & $2^n$ & $\F_2^n$ & $1$ & $1$ \\
 \hline
\end{tabular}   \vspace*{3mm}
\caption{Rules with trivial dynamics group.\label{tbl:trivial}}
\end{table}

\begin{defn}[Avoiding words]
  For each $n$, let $\no_{XY\ldots}$ denote the set of states in
  $\F_2^n$ that do not contain any subwords of the form $X$, $Y$,
  \ldots.  (The letter $\no$ stands for ``no''.) For example
  $\no_{\mbox{`11'}}$ is the collection of states without adjacent
  $1$s, keeping in mind that we view the subscripts mod $n$.
\end{defn}

\begin{defn}[Abbreviations]
  To simplify notation, we introduce six abbreviations: $A$=`11',
  $B$=`000', $C$=`111', $D$=`010', $E$=`101', and $F$=`1100'.  Thus,
  $\no_{AE}$ represents the states in $\F_2^n$ with no subwords of the
  form `11' or `101'.  In addition, let $\zer$ and $\one$ refer to the
  state with all $0$s and all $1$s, respectively.
\end{defn}

\begin{rem}[Why these words]
  The words we have chosen to abbreviate are those needed to
  efficiently describe the periodic sets of the $41$ representative
  rules.  The words $B$, $C$, $D$ and $E$ are triples that need to be
  avoided, while the words $A$ and $F$ deserve additional explanation.
  Avoiding the triple `110' is equivalent to avoiding the subword `11'
  while allowing the state $\one$, and avoiding the triple `011' leads
  to the same conditions.  It thus makes sense to abbreviate the word
  `11' and treat the state $\one$ separately.  The word $F$ = `1100'
  is only needed to describe $P_{28}$ and $P_{29}$, so we postpone our
  discussion of this abbreviation until Section~\ref{sec:exceptional}.
  We note that this is the only abbreviated word that is not
  left-right symmetric.
\end{rem}

The periodic sets for the $26$ representative Wolfram rules under
discussion, calculated as described above, are listed in
Table~\ref{tbl:trivial}.  The tag of rule $k$, along with the decimal
of its inversion, reflection and inversion-reflection are also
included.  Finally we turn to a calculation of the number of orbits
for each of these rules.

\begin{rem}[Recurrence relations]
  Because the dynamics groups are trivial, we have $o_k = \card{P_k}$
  in each case.  Moreover, since each of our sets is defined by a
  finite list of configurations that it avoids, it is well-known that
  the number $a_n$ of acceptable configurations for each $n$ are the
  coefficients of an easily calculated rational generating
  function~\cite{Stanley:97}.  As a consequence, they satisfy a
  constant coefficient recurrence relation. The last statement also
  follows from~\cite[Theorem~5.3, p.~132]{Mortveit:07}.
\end{rem}

\begin{table}
\small
  \begin{tabular}{|c|c|c|l|c|}\hline
 name & words to avoid & tiles & recurrence relation & Sloane \\  \hline\hline
  $\no_A$ & `11' & `0', `10' & $a_n=a_{n-1}+a_{n-2}$ & A000032 \\ \hline
  $\no_{AB}$ & `11', `000' & `10', `100' & $a_n=a_{n-2}+a_{n-3}$ & A001608 \\ \hline
  $\no_{AE}$ & `11', `101' & `0', `100' & $a_n=a_{n-1}+a_{n-3}$ & A001609 \\  \hline
  $\no_{BCF}$ & `000', `111', `1100' & `10', `100', `110' & $a_n=a_{n-2}+2a_{n-3}$ & A072328 \\ \hline
  $\no_{BC}$ & `000', `111' & `10', `100', `110', `1100' & $a_n=a_{n-2}+2a_{n-3}+a_{n-4}$ & A007040 \\ \hline
  $\no_{C}$ & `111' & `0', `10', `110' & $a_n=a_{n-1}+a_{n-2}+a_{n-3}$ & A001644 \\ \hline
  $\no_{CD}$ & `111', `010' & `0', `110' & $a_n=a_{n-1}+a_{n-3}$ & A001609 \\ \hline
  $\no_{D}$ & `010' & --- & $a_n=2a_{n-1}-a_{n-2}+a_{n-3}$ & A109377 \\ \hline
  $\no_{DE}$ & `010', `101' & --- & $a_n=2a_{n-1}-a_{n-2}+a_{n-4}$ & A007039 \\ \hline
  \end{tabular}  \vspace*{3mm}
  \caption{Recurrence relations for the basic periodic sets.}
  \label{tbl:recurrence}
\end{table}

In all but two cases these negative descriptions can be reformulated
as positive ones that allow us to compute the recurrence with ease.
For example, the states in $\no_A$ are those without adjacent $1$s and
every such state can be uniquely decomposed into subwords of the form
`0' and `10' that we call ``tiles''.  Counting these involve the Lucas
numbers and, indirectly, the Fibonacci numbers.

\begin{defn}[Fibonacci and Lucas]\label{def:fib-luc}
  In order to count the size of $\no_A$ we first consider the number
  $b_n$ of ways to build a word of length $n$ out of the tiles $0$ and
  $10$, i.e.  with no cyclic subscripting.  By focusing on the type of
  the final tile, we see that the number of such tiles satisfy the
  recurrence relation $b_n = b_{n-1} + b_{n-2}$ with initial
  conditions $b_1=1$ and $b_2=2$.  The unique solution of the
  recurrence is $b_n=\fib_{n+1}$ where $\fib_n$ are the famous
  \emph{Fibonacci numbers} with values $\{1,1,2,3,5,8,13,21,\ldots\}$
  starting with $\fib_1$.  Returning to the cyclic version, there are
  three ways the vertex $i$ can be covered by a tile: it can be a $0$
  tile, the first digit of a $10$ tile or the second digit of a $10$
  tile.  Once the tile containing vertex $i$ has been placed, the
  remaining problem involves tiling a word.  Thus $a_n = b_{n-1} +
  2b_{n-2}$.  It is now easy to see that $a_n$ satisfies the same
  recurrence relation as $b_n$ but with different initial conditions.
  In this case $a_1=1$ and $a_2=3$ and the solution of the recurrence
  is $a_n = \luc_n$ where $\luc_n$ are the nearly as famous
  \emph{Lucas numbers} with values $\{1,3,4,7,11,18,29,\ldots\}$
  starting with $\luc_1$.
\end{defn}

This analysis of the size of $\no_A$ easily extends to those cases
with a tiling description supplementing the description by words to
avoid.  The tiles in each case and the corresponding recurrence
relation are listed in Table~\ref{tbl:recurrence}.  The final column
is a reference to the appropriate entry in Neil Sloane's Online
Encyclopedia of Integer Sequences~\cite{Sloane:06}.  The remaining two
cases, $\no_D$ and $\no_{DE}$, are classical situations where isolated
$1$s and/or isolated $0$s are to be avoided.  See the references
listed in their entries in \cite{Sloane:06} for further details.

%%%%%%%%%%%%%%%%%%%%%%%%%%%%%%%%%%%%%%%%%%%%%%%%
\section{Invertible rules}\label{sec:invertible}
%%%%%%%%%%%%%%%%%%%%%%%%%%%%%%%%%%%%%%%%%%%%%%%%

We now turn our attention to Wolfram rules where every state is
periodic. Of the $41$ representative Wolfram rules, $9$ of them fall
into this category.  Our results for these $9$ rules are summarized in
Table~\ref{tbl:invertible}.  Note that the more complicated numbers
and groups $o_k$ and $G_k$ are not explicitly listed in the table and
only described in the text.  As a convention, when it is clear that
rule $k$ is the rule under discussion, we use $F_i$ instead of the
more cumbersome $\Wolf_i^{(k)}$ to denote the rule that updates the
$i^{th}$ coordinate.  Finally, recall that the Coxeter label $c_k$
must be a divisor of $12$ and note that all six possible values occur
among these ``invertible'' rules.

\begin{defn}[Invertible rules]\label{def:invertible}
  An \sds $\Fy$ is called \emph{invertible} when every state $\F^n$ is
  periodic, or equivalently when every update rule $F_i$ is a
  bijection.  To see the equivalence note that when every update rule
  is bijective, their composition is a permutation and every state is
  periodic.  Conversely, if every state is periodic, then the \sds
  maps are permutations and the only way this can happen is when every
  update rule is at least injective and hence bijective.
\end{defn}

It is also worth noting that when there are only two possible vertex
states, invertibility is characterized by the absence of the symbols
{\tt 0}s and {\tt 1}s in the tags of the update rules.  Before
discussing the invertible rules individually, we need some additional
notation.

\begin{defn}[Blocks]
  If $\y$ is any state other than $\zer$ or $\one$ then it consists of
  alternating strings of $0$s an $1$s that we abbreviate using
  exponents.  For example $1^50^3$ is shorthand for the word
  $11111000$.  We call a maximal subword of the form $0^i$, $1^j$,
  $0^i1^j$ or $1^j0^i$ (with $i$ and $j$ positive) a \emph{$0$-block},
  a \emph{$1$-block}, a \emph{$01$-block} and a \emph{$10$-block},
  respectively.  For example, keeping in mind the cyclic nature of the
  subscripts, the state $\y = 1101001$ has two $1$-blocks ($1$ and
  $1^3$), two $0$-blocks ($0$ and $0^2$), two $01$-blocks ($01$ and
  $0^21^3$), and two $10$-blocks ($10^2$ and $1^30$).  Every block has
  a length. Blocks of length~$1$ are called \emph{isolated} and longer
  blocks are called \emph{nontrivial}.
\end{defn}

The effect of each of the four possible {\tt x}'s in a tag can be
described in this language.

\begin{rem}[Tags and dynamics]%%their effects]
  When $t_0$={\tt x} ($010 \leftrightarrow 000$) isolated $1$s can be
  created and removed and when $t_3$={\tt x} ($111 \leftrightarrow
  101$) isolated $0$s can be created and removed.  When $t_1$= {\tt x}
  ($011 \leftrightarrow 001$) the boundary between a $1$-block and the
  $0$-block to its left can be shifted left or right and when $t_2$=
  {\tt x} ($110 \leftrightarrow 100$) the boundary between a $1$-block
  and the $0$-block to its right can be shifted left or right.
\end{rem}

\begin{table}
\begin{tabular}{|c|c|rrr|c|c|c|c|}
  \hline
  rule & tag & inv & refl& i.r. & $o_k$ & $P_k$ & $c_k$ & $G_k$  \\
  \hline \hline
  204 & \tt{----} & 204 & 204 & 204 & $2^n$ & $\F_2^n$ & $1$ & $1$ \\ \hline
  51 & \tt{xxxx} &  51 &  51 &  51 &  $1$ & $\F_2^n$ & $2$ & $\Z_2^n$ \\ \hline
  60 & \tt{xx--} & 195 & 102 & 153 &  $2$ & $\F_2^n$ &  $4$ & $\psl_n(\Z_2)$ \\ \hline
  150 & \tt{-xx-} & 150 & 150 & 150 & $\floor{\frac{n}{2}}+2$ & $\F_2^n$ & $3$ & Thm~\ref{thm:150} \\
  105 & \tt{x--x} & 105 & 105 & 105 & $o_{105}$ & $\F_2^n$ & $3$ & Thm~\ref{thm:105} \\ \hline
  156 & \tt{-x--} & 198 & 198 & 156 & $\card{\no_A\cup\one}$ & $\F_2^n$ & $6$ & Thm~\ref{thm:156} \\
  201 & \tt{---x} & 108 & 201 & 108 & $\card{\no_D}$ & $\F_2^n$ & $6$ & Thm~\ref{thm:201} \\ \hline
  57 & \tt{xx-x} &  99 &  99 &  57 &  $1$ & $\F_2^n$ & $12$ & Conj~\ref{conj:57} \\
  54 & \tt{xxx-} & 147 &  54 & 147 &  $2$ & $\F_2^n$ & $12$ & Conj~\ref{conj:54} \\ \hline
\end{tabular} \vspace*{2mm}
\caption{Invertible rules.\label{tbl:invertible}}
\end{table}

We now discuss the invertible rules one at a time.

%%%%%%%%%%%%%%%%%%%%%%%%
\subsection*{Rule $204$}
As we already noted in the last section, under rule $204$, with tag
{\tt ----}, all of the local functions leave the state of the system
unchanged, every local function induces the trivial permutation, and
the group generated is the trivial group.  There are thus $2^n$ orbits
since every distinct state is an orbit.

%%%%%%%%%%%%%%%%%%%%%%%
\subsection*{Rule $51$}

Under rule $51$, with tag {\tt xxxx}, the local functions ignore their
context and always alter the value.  Since the local functions
pairwise commute, $c_{51}=2$, and the dynamics group is a quotient of
$\Z_2^n$.  But since the composition of every distinct subset of local
functions toggles a distinct subset of vertex states, there are at
least $2^n$ elements in $G_{51}$.  Thus $G_{51}$ is isomorphic to
$\Z_2^n$.  Finally, it is easy to see that there is only one orbit, so
$o_{51}=1$.

%%%%%%%%%%%%%%%%%%%%%%
\subsection*{Rule $60$}

The dynamics group $G_{60}$ is interesting because its structure is
slightly unexpected.  The key observation is that when vertex $i$ is
updated, its new value is its old value plus the value of the vertex
immediately to its left.  In other words $\Wolf_i^{(60)}$ replaces
$y_i$ with $y_i+y_{i-1}$, which leads to a matrix representation of
the update rules.  Viewing $\y$ as a column vector, the effect of
updating vertex $i$ can be achieved by multiplying $\y$ on the left by
the matrix $\A_i:=I + E_{i,i-1}$ where $I$ is the $n\times n$ identity
matrix and $E_{i,j}$ is the elementary matrix with $0$s everywhere
except for a single $1$ in the $(i,j)$ position.  Matrix
multiplication by $\A_i$ is a concise description of the function
$\Wolf_i^{(60)}$ from $\F_2^n$ to $\F_2^n$ and thus the matrix group
generated by the $\A_i$s is isomorphic to the dynamics group $G_{60}$.
Since each $\A_i$ has determinant $1$, it is clear that $G_{60}$ is a
subgroup of $\psl_n(\Z_2)$, and, in fact, it is well-known that these
matrices generate all of $\psl_n(\Z_2)$~\cite[pg. 455]{Weir:55}.  Thus
$G_{60}$ is isomorphic to $\psl_n(\Z_2)$.  Under this group action it
is clear that $\zer$ is fixed and the remaining states form a single
orbit.  A calculation shows $c_{60} = 4$.

%%%%%%%%%%%%%%%%%%%%%%%%
\subsection*{Rule $150$}

Rule $150$ is similar to rule $60$ but this time, when vertex $i$ is
updated, its new value is the sum of its old value plus the value of
the vertices immediately to its left and its right.  This is better
known as the \emph{parity function}.  More explicitly, $F_i =
\Wolf_i^{(150)}$ replaces $y_i$ with $y_{i-1}+y_i+y_{i+1} \mod 2$,
which leads to a matrix representation of the update rules as before.
If we define the matrices $\A_i$ as $I + E_{i,i-1} + E_{i,i+1}$ then
the group $G_{150}$ can be identified as the subgroup of
$\psl_n(\Z_2)$ the matrices $\A_i$ generate.  Unlike rule $60$ it is
not clear which subgroup of $\psl_n(\Z_2)$ this generates.  An
approach via Coxeter groups is more successful.

Since an easy calculation shows that $c_k=3$, $G_{150}$ is a quotient
of $\cox_n$.  Recall from Remark~\ref{rem:cox} that $\cox_n \cong \Zo
\rtimes \sym_n$ as well as the conventions established there.  The
analysis of $G_{150}$ involves two steps.  The first is to show that
the map $\cox_n\to G_{150}$ is a factor of the map $\cox_n \to \sym_n$
(i.e. the later map decomposes as $\cox_n\to G_{150} \to \sym_n$).  To
see that $G_{150}$ maps onto $\sym_n$ in a manner consistent with the
projection $\cox_n\to \sym_n$ start with a state $\y$ and imagine that
the numbers $1$ up to $n$ are placed in the gaps between the $n$
positions.  In particular, initially place the number $i$ between
$y_{i-1}$ and $y_i$.  When $F_i$ is applied to a state $\y$ switch the
numbers on either side of $y_i$ fixing all the others, in addition to
updating the value of $y_i$.  We claim that if $i$ was originally in a
gap that marked the end of a $0$-block or $1$-block, then the same is
true of the place where $i$ ends up in the final state.  This follows
easily from the way rule $150$ updates states.  From this it is not
too hard to see that if a sequence of update rules fixes every state,
then each of the numbers $1$ through $n$ must also return to their
original position.  More concretely, if the numbers $1$ through $n$ do
not all return to their original positions, it is easy to find an
explicit state that is not fixed by this sequence of update rules.
This means that there is a well-defined group homomorphism $G_{150}$
to $\sym_n$ that sends the permutation $[\Fy,\pi] = F_{\pi_n} \circ
\cdots \circ F_{\pi_1}$ to its permutation of the set
$\{1,2,\ldots,n\}$.  Since this matches the image of
$s_{\pi_n}\cdots\circ s_{\pi_1}$ under the map $\cox_n\to \sym_n$ we
have the factorization we desire.

The existence of maps $\cox_n \to G_{150} \to \sym_n$ imply that the
kernel of the first map consists solely of pure translations and we
only need to analyze which pure translations lie in the kernel in
order to completely understand the group $G_{150}$.  To do this we use
the concrete description of $\cox_n$ given in Remark~\ref{rem:cox}.
One of the generating pure translations in $\cox_n$ is the element
$T_1= s_1 s_2 \cdots s_{n-1} s_n s_{n-1} \cdots s_3 s_2$.  The image
of this inside $G_{150}$ is the element $F_1 \circ F_2 \circ \cdots
\circ F_{n-1} \circ F_n \circ F_{n-1} \cdots \circ F_3 \circ F_2$.  If
we apply this sequence of update rules to an arbitrary state $\y$
(using the parity function as we should), the final result is $\y +
(y_2+y_n)\one$.  More generally, the image of $T_i$ acts on states by
sending $\y$ to $\y + (y_{i-1} + y_{i+1})\one$.  Applying $T_i$ twice
is clearly trivial so the kernel contains the subgroup of $\Zo$ that
the vectors $2T_i$ generate.  This set is $(2\Z)^n \perp \one$.  As a
consequence, $G_{150}$ is a quotient of the group $\Zt \rtimes \sym_n$
of size $2^{n-1} \cdot n!$.  When $n$ is even there is another pure
translation in the kernel, namely, the result of applying once each
$T_i$ with an odd subscript.  As a pure translation this adds the
vector $\langle -1,1,-1,\cdots,-1,1\rangle$.  Equivalently, when $n$
is even $G_{150}$ is a quotient of $\left(\Zt/\langle \one \rangle
\right) \rtimes \sym_n$.  (Note that adding $\one$ to a state $\x$
replaces $\x$ with its \emph{complement}, the $0$s become $1$s and
vice versa.  Moreover, this pure translation commutes with the
symmetric group action and is central in $\Zt \rtimes \sym_n$.)  We
now show that these are the only pure translations in the kernel.

\begin{thm}[Rule $150$]\label{thm:150}
  When $n$ is odd the group $G_{150}$ is isomorphic to $\Zt \rtimes
  \sym_n$ and when $n$ is even it is isomorphic to $\left(\Zt/\langle
    \one \rangle \right) \rtimes \sym_n$.  In particular,
  $\card{G_{150}} = 2^{n-1}\cdot n!$ when $2\nmid n$ and $2^{n-2}\cdot
  n!$ when $2\mid n$.
\end{thm}

\begin{proof}
  That $G_{150}$ is a quotient of these groups was shown above, so we
  only need to show that we have found the full kernel of the map
  $\cox_n\to G_{150}$.  To see whether there are any other pure
  translations that are trivial in $G_{150}$ note that $T_{i-1}$ and
  $T_{i+1}$ are the only generating pure translations that add $y_i
  \one$ to $\y$.  Thus, if $a_1 T_1 + a_2 T_2 + \cdots + a_n T_n$ acts
  trivially on every state $\y$, the parity of $a_{i-1}$ and $a_{i+1}$
  must match for every $i$.  When $n$ is odd, this means that all the
  $a_i$ have the same parity and every such element is one we already
  know lies in the kernel.  When $n$ is even there is one additional
  possibility.  Perhaps the $a_i$s with even subscripts have one
  parity and the ones with odd subscripts have the other.  Removing
  summands we already know to lie in the kernel and using the relation
  $T_1 + T_2 +\cdots + T_n = \zer$ if necessary, we see that this
  possibility is equivalent to $T_1 + T_3 + T_5 + \cdots + T_{n-1}$,
  again, an element we already know lies in the kernel.
\end{proof}

Although we did not need its orbit structure in order to analyze
$G_{150}$, it is easy to see that the changes rule $150$ allows (and
the only changes it allows) are the alterations of the boundaries of
the $0$-blocks and $1$-blocks.  Thus two states belong to the same
orbit iff they have the same number of $0$-blocks and the same number
of $1$-blocks.  In particular $o_{150} = \floor{\frac{n}{2}}+2$ with
the $2$ corresponding to the fixed states $\zer$ and $\one$.

%%%%%%%%%%%%%%%%%%%%%%%%
\subsection*{Rule $105$}

Rule $105$ is the negation of the parity function and its analysis is
very similar to our analysis of rule $150$.  In particular,
$F_i=\Wolf_i^{(105)}$ replaces $y_i$ with $1+y_{i-1}+y_i+y_{i+1} \mod
2$.  The group $G_{105}$ can be described as a subgroup of a linear
group (with an extra row and column added for constants) but viewing
it as a Coxeter quotient is more fruitful.  The Coxeter label
$c_{105}=3$, and $\cox_n \to G_{105} \to \sym_n$ as above.  The
argument for such a factorization is similar in spirit to the one
presented above for rule $150$ but complicated by the presence of
negations.  We omit the details.

The existence of maps $\cox_n \to G_{105} \to \sym_n$ imply that the
kernel of the first map consists solely of pure translations and we
only need to analyze which pure translations lie in the kernel in
order to completely understand the group $G_{105}$.  This time the
image of the pure translation $T_i$ in $G_{105}$ acts on states by
sending $\y$ to $\y + (1 + y_{i-1} + y_{i+1})\one$.  Applying $T_i$
twice is clearly trivial so the kernel contains the subgroup of $\Zo$
that the vectors $2T_i$ generate.  This set is $(2\Z)^n \perp \one$.
As a consequence, $G_{105}$ is a quotient of the group $\Zt \rtimes
\sym_n$ of size $2^{n-1} \cdot n!$.  When $n$ is a multiple of $4$
there is another pure translation in the kernel, namely, the result of
applying once each $T_i$ with an odd subscript.  The reason $n$ needs
to be a multiple of $4$ and not merely even is that when $n$ is twice
an odd number $\y$ is sent to $\y + \one$.  As a pure translation this
adds the vector $\langle -1,1,-1,\cdots,-1,1\rangle$.  Equivalently,
when $n$ is even $G_{105}$ is a quotient of $\left(\Zt/\langle \one
  \rangle \right) \rtimes \sym_n$.  We now show that these are the
only pure translations in the kernel.

\begin{thm}[Rule $105$]\label{thm:105}
  When $n$ is not a multiple of $4$ the group $G_{105}$ is isomorphic
  to $\Zt \rtimes \sym_n$ and when it is a multiple of $4$ it is
  isomorphic to $\left(\Zt/\langle \one \rangle \right) \rtimes
  \sym_n$.  In particular, $\card{G_{105}} = 2^{n-1}\cdot n!$ when
  $4\nmid n$ and $\card{G_{105}} = 2^{n-2}\cdot n!$ when $4\mid n$.
\end{thm}

\begin{proof}
  That $G_{105}$ is a quotient of these groups was shown above, so we
  only need to show that we have found the full kernel of the map
  $\cox_n\to G_{105}$.  To see whether there are any other pure
  translations that are trivial in $G_{105}$ note that $T_{i-1}$ and
  $T_{i+1}$ are the only generating pure translations that add $y_i
  \one$ to $\y$.  Thus, if $a_1 T_1 + a_2 T_2 + \cdots + a_n T_n$ acts
  trivially on every state $\y$, the parity of $a_{i-1}$ and $a_{i+1}$
  must match for every $i$.  When $n$ is odd, this mean that all the
  $a_i$ have the same parity and every such element is one we already
  know lies in the kernel.  When $n$ is even there is one additional
  possibility.  Perhaps the $a_i$s with even subscripts have one
  parity and the ones with odd subscripts have the other.  Removing
  summands we already know to lie in the kernel and using the relation
  $T_1 + T_2 +\cdots + T_n = \zer$ if necessary, we see that this
  possibility is equivalent to $T_1 + T_3 + T_5 + \cdots + T_{n-1}$.  If
  $n$ is twice an odd number then $\y$ is sent to $\y + \one$ and this
  is not in the kernel, but when $n$ is twice an even number (i.e. $4
  \mid n$) then this is an element we already know lies in the kernel.
\end{proof}

Finally, for the sake of completeness, we include a description of
$o_{105}$ but leave a verification of our assertions as an exercise.
When $n$ is twice an odd number $o_{105}=2 \floor{\frac{n}{4}}+4$,
when $n$ is twice an even number $o_{105}=2 \floor{\frac{n}{4}}+2$,
and when $n$ is odd, $o_{105}=2$ (and these are the orbits of $\zer$
and $\one$).

%%%%%%%%%%%%%%%%%%%%%%%%
\subsection*{Rule $156$}

Wolfram rule $156$ is a case where the orbit structure is useful and
we use it to prove that $G_{156}$ is a subdirect product of symmetric
groups of particular sizes in a very concrete fashion.  The first
thing to notice is that the only changes rule $156$ allows are growing
or shrinking $1$-blocks from the right.  Isolated $1$s and $0$s can
neither be created nor removed (so the number of blocks is an
invariant) and the left end of a $1$-block is never moved (so that
subwords of the form $01$ persist forever).  To analyze the number of
orbits, apply the rewrite rule $110\rightarrow 100$ to remove as many
$1$s as possible.  It should be clear that (unless we started with the
state $\one$) this process ends when all remaining $1$s are isolated.
Moreover, these remaining $1$s are precisely the leftmost $1$s in the
initial blocks so there is a \emph{unique} state in each orbit with
isolated ones.  In other words, the orbits are in natural bijection
with the set $\no_A \cup \one$, and thus $o_{156} = \card{\no_A \cup
  \one} = \luc_n+1$.

We now examine the way the updates rules act on the states in each
orbit individually.  The fixed state $\one$ can be discarded since the
trivial group it generates does not contribute meaningfully to the
subdirect product.  Next we note that the dynamics group restricted to
a single orbit usually splits further as a direct product.  For
example, consider the state $\y = 1000001000$. It has a two
$10$-blocks, one of length $6$ and one of length $4$.  The orbit of
$\y$ consists of all words that can be written as a $10$-block of
length $6$ followed by a $10$-block of length $4$.  In particular, the
first six digits can be $100000$, $110000$, $111000$, $111100$ or
$111110$ and the last four can be $1000$, $1100$ or $1110$.  These
fifteen combinations form the complete orbit.  In addition, writing
$F_i$ in place of $\Wolf_i^{(156)}$, the only update rules that act
nontrivially on this orbit are $F_2$, $F_3$, $F_4$, $F_5$, $F_8$ and
$F_9$, but $F_2$, $F_3$, $F_4$ and $F_5$ commute with $F_8$ and $F_9$
so $G_{156}$ restricted to this orbit splits as a direct product
$\langle F_2, F_3, F_4, F_5 \rangle \times \langle F_8, F_9\rangle$.

If we focus in on the first six digits for a moment and name the
possibilities $x_1 = 100000$, $x_2 = 110000$, $x_3=111000$,
$x_4=111100$ and $x_5=111110$ then we see that $F_2$ swaps $x_1$ and
$x_2$ and fixes the other $x_i$, $F_3$ swaps $x_2$ and $x_3$, $F_4$
swaps $x_3$ and $x_4$ and $F_5$ swaps $x_4$ and $x_5$.  The induced
permutations of the five possibilities are $F_2 = (1 2)$, $F_3 = (2
3)$, $F_4 = (3 4)$ and $F_5 = (4 5)$.  Thus, these four update rules,
restricted in this way generate a copy of $\sym_5$ with the functions
$F_2$, $F_3$, $F_4$ and $F_5$ acting as the standard Coxeter
generating set.

More generally, if a state $\y$ contains a $10$-block of length $m+1$
starting at position $\ell$, then every state in the orbit of $\y$
contains a $10$-block of this length at this location. There are $m$
possibilities for a block of this type that we call $x_i = 1^i0^j$
where $i,j>0$ and $i+j=m+1$.  The update rules $F_{\ell+1}$, \ldots,
$F_{\ell+m-1}$ are the only ones that alter this block and they
generate a copy of $\sym_m$ acting on the possibilities $x_i$ in the
standard way.  We call this group $\sym_m^{(\ell)}$.  As a group it is
$\sym_m$ and the number $\ell$ indicates which update rules generate
$\sym_m$ and in what way.  Returning to our earlier example, we see
that the quotient of $G_{156}$ obtained by restricting to the orbit of
$\y = 1000001000$ is the group $\sym_5^{(1)} \times \sym_3^{(7)}$.
More generally, if $\y$ has $10$-blocks of lengths $m_1+1$, $m_2+1$,
\ldots, $m_k+1$ starting at $\ell_1$, $\ell_2$, \ldots, $\ell_k$ then
the quotient of $G_{156}$ restricted to the orbit of $\y$ is the group
$\sym_{m_1}^{(\ell_1)} \times\sym_{m_2}^{(\ell_2)} \times \cdots
\times \sym_{m_k}^{(\ell_k)}$.

To summarize, $G_{156}$ is a subdirect product of groups, each of which is
isomorphic to a direct product of groups of the form
$\sym_m^{(\ell)}$.  When written out completely using the orbit
structure, the groups $\sym_m^{(\ell)}$ occur multiple times and the
projection of $G_{156}$ onto each of these repeated factors is
identical.  As a consequence, the group $G_{156}$ can be embedded in a
direct product where each of these groups only occurs once.  Note that the
position $\ell$ can be any number $1$ through $n$, but that $m+1$ can
only be $2$ through $n$ excluding $n-1$ since the complement of the
$10$-block with length $m+1$ must be tilable by $10$-blocks and these
have length at least $2$.

\begin{thm}[Rule $156$]\label{thm:156}
  The group $G_{156}$ is a subdirect product of symmetric groups.  In
  particular, \[ G_{156} \subset \prod_{\ell=1}^{n} \left(
    \sym_{n-1}^{(\ell)} \times \prod_{m=2}^{n-3} \sym_m^{(\ell)}
  \right) \] and each generator of $G_{156}$ restricted to a factor is
  either trivial or a standard Coxeter generator as described above.

\end{thm}

Note that the factors with $m=1$ have been eliminated from these
products, which is possible because the groups $\sym_1^{(\ell)}$ are
trivial.  We computed the order of $G_{156}$ for $4\leq n\leq 7$, and
the results are shown in Table~\ref{tbl:73-201-156}.  Notice that for
$n>4$ each of these orders are $n^{th}$ powers.

\begin{table}
\centering
\begin{tabular}{|l|c|c|c|} \hline
  $n$ & {\small $|G_{156}|$} & {\small $|G_{201}|$} & {\small $|G_{73}|$} \\ \hline
  $4$ & $2^3.3^4$ &  $(7!/2)$ & $(7!/2)$ \\
  $5$ & $2^{15}.3^5$ & $2^5.(11!/2)$ & $2^5.(11!/2)$ \\
  $6$ & $2^{18}.3^{12}.5^6$ & $2^6.3^6.(18!/2)$ & $2^5.3^6.(18!/2)$ \\
  $7$ & $2^{42}.3^{28}.5^7$ & $2^{21}.3^{14}.5^7.(29!/2)$  &  $2^{21}.3^7.5^7.(29!/2)$\\ \hline
\end{tabular} \vspace*{2mm}
\caption{The orders of $G_{156}$, $G_{201}$ and $G_{73}$ for small
  values of $n$.\label{tbl:73-201-156}}
\end{table}

%%%%%%%%%%%%%%%%%%%%%%%%
\subsection*{Rule $201$}

Rule $201$ is a second situation where the orbit structure can be used
to simplify our analysis.  In particular we use orbit structures to
prove that $G_{201}$ is a subdirect product of groups that are all
(conjecturally) symmetric or alternating groups of particular sizes in
a very concrete fashion.  The first step is to analyze the orbit
structure.  The only changes under rule $201$ are the creation or
removal of isolated $1$s.  In particular the size and location of
nontrivial $1$-blocks in a state $\y$ are an invariant of its orbit.
In fact, since we can remove all isolated $1$s, any two states with
the same set of nontrivial $1$-blocks belong to the same orbit and
every orbit has exactly one state with no isolated $1$s.  Thus
$o_{201} = \card{\no_D}$, the number of states with no isolated $1$s.

Let $\y$ be a state with no isolated $1$s.  It should be clear that
the restriction of the group $G_{201}$ to the orbit of $\y$ will split
as a direct product with one factor for each $0$-block in $\y$.  This
is because all of the local functions that insert and remove isolated
$1$s in one $0$-block will commute with those that insert and remove
isolated $1$s in a different $0$-block, being separated by nontrivial
$1$-blocks.  As was the case with rule $156$, the orbit groups of
$G_{201}$ split into factors, the refined direct product contains
redundancies and these redundancies can be removed to give a better
subdirect product representation of $G_{201}$.  In particular,
$G_{201}$ embeds in the subdirect product of orbit groups acting on
orbits with only one $0$-block.  Let $\lucgp_n$ denote the group
obtained by restricting $G_{201}$ to the orbit of $\zer$ and let
$\fibgp_m^{(\ell)}$ be the group obtained by restricting $G_{201}$ to
the orbit of the state $\y$ with a single $0$-block of length $m$
starting at position $\ell$.  The notations stand for \emph{Lucas
  group} and \emph{Fibonacci group} and are suggested by the fact that
the orbit of $\zer$ under $G_{201}$ is $\no_A$ of size $\luc_n$ and
the orbit of $\y$ under $G_{201}$ has size $\fib_m$ (the first $0$ is
fixed and the remaining $m-1$ digits of the $0$-block are tiled by
tiles $0$ and $10$).  In the case of $\lucgp_n$ we can be slightly
more precise.

\begin{prop}[Lucas groups]\label{prop:lg}
  The group $\lucgp_n$ is a subgroup of $\sym_{\luc_n}$ when
  $\fib_{n-1}$ is odd and a subgroup of $\alt_{\luc_n}$ when
  $\fib_{n-1}$ is even.
\end{prop}

\begin{proof}
  Since $\lucgp_n$ acts on a set of size $\card{\no_A}=\luc_n$, only
  the second assertion needs to be established.  Let $F_i^*$ denote
  the restriction of $\Wolf_i^{(201)}$ to $\no_A$ and note that these
  permutations generate $\lucgp_n$.  To count how many $2$-cycles are
  in the disjoint cycle notation for $F_i^*$, note that for each
  $\y\in \no_A$, $F_i^*(\y)\neq\y$ iff $y_{i-1}=y_{i+1}=0$.  Thus, the
  number of $2$-cycles in $F_i^*$ equals the number of words $y_{i+1}
  y_{i+2} \ldots y_{i-2} y_{i-1}$ of length $n-1$ with only isolated
  $1$s and $y_{i+1}=y_{i-1}=0$, which equals the number of ways to
  tile a word $y_{i+2} \ldots y_{i-2} y_{i-1}$ of length $n-2$ with
  tiles $0$ and $10$.  As we saw in Definition~\ref{def:fib-luc}, this
  number is $\fib_{n-1}$.  Thus $F_i^*$ is an even permutation iff
  $\fib_{n-1}$ is even, and when this is true, $\lucgp_n$ is a
  subgroup of $\alt_{\luc(n)}$.
\end{proof}

The following conjecture is based on computational evidence for small
values of $n$.  We have checked, for example, that $\fibgp_6^{(\ell)}
\cong \sym_8 = \sym_{\fib_6}$ and that $\lucgp_4 \cong \alt_7 =
\alt_{\luc_4}$ (as predicted since $\fib_3 = 2$ is even).

\begin{conj}[Fibonacci groups and Lucas groups]\label{conj:fg-lg}
  For all $m$, $\fibgp_m^{(\ell)} \cong \sym_{\fib_m}$ and for all
  $n$, $\lucgp_n \cong \sym_{\luc_n}$ when $\fib_{n-1}$ is odd and
  $\lucgp_n \cong \alt_{\luc_n}$ when $\fib_{n-1}$ is even.
\end{conj}

The parity of the Fibonacci numbers is quite predictable and a more
direct statement is possible, but the one given is more closely tied
to the reason $\lucgp_n$ stays inside $\alt_{\luc_n}$.  The
obstruction to establishing Conjecture~\ref{conj:fg-lg} is how far the
permutations generating the Fibonacci groups and Lucas groups differ
from the standard generating sets of the symmetric and alternating
groups.  In any case, this analysis of the orbit structure under rule
$201$ gives the following subdirect product decomposition for
$G_{201}$.  The values of $m$ only range from $3$ to $n-2$ to exclude
trivial groups on the low end and because the complement of the
$0$-block needs to leave room for a nontrivial block of $1$s.

\begin{thm}[Rule $201$]\label{thm:201}
  The group $G_{201}$ is a subdirect product of Lucas groups and
  Fibonacci groups.  In particular, \[ G_{201} \subset \lucgp_n \times
  \prod_{\ell=1}^{n} \prod_{m=3}^{n-2} \fibgp_m^{(\ell)} \]
\end{thm}

Finally, we used a computer program to calculate the order of
$G_{201}$ for $4\leq n\leq 7$, and the results are shown in
Table~\ref{tbl:73-201-156}.  Notice that the answer in each case is an
$n^{th}$ power times the size of $\alt_{\luc_n}$.

%%%%%%%%%%%%%%%%%%%%%%%
\subsection*{Rule $57$}

Wolfram rule $57$ can introduce and remove isolated $0$s and $1$s and
it can grow and shrink $0$-blocks from the left and $1$-blocks from
the right.  With so much flexibility it is easy to see that there is
only one orbit and thus $o_{57}=1$.  Moreover, because
$\Wolf_i^{(57)}(\y) = \y$ only when $y_{i-1}=0$ and $y_{i+1}=1$, the
number of $2$-cycles in the disjoint cycle representation of
$\Wolf_i^{(57)}$ is $2^{n-3}$.  In particular, it is an even
permutation for $n>4$ and we conclude that $G_{57}$ is a subgroup of
$\alt_{2^n}$.

\begin{conj}[Rule $57$]\label{conj:57}
  For all $n>4$, $G_{57} \cong \alt_{2^n}$.
\end{conj}

We have verified Conjecture~\ref{conj:57} explicitly up to $n=8$.  The
reason Conjecture~\ref{conj:57} is difficult to establish abstractly
is that the permutations $\Wolf_i^{(57)}$ are far removed from the
standard generating sets of $\alt_{2^n}$.

%%%%%%%%%%%%%%%%%%%%%%%
\subsection*{Rule $54$}

Rule $54$ is similar to rule $57$ except that the state $\zer$ is now
fixed. Since $1$-blocks can grow to the left or the right and isolated
$0$s can be removed, all states other than $\zer$ are in the same
orbit as $\one$.  Thus $o_{54}=2$.  Since fixed states only contribute
trivial groups to the subdirect product structure, they can be ignored
when computing the dynamics group.  In particular, $G_{54}$ is a
subgroup of $\sym_{2^n-1}$.  Because $\Wolf_i^{(54)}(\y) = \y$ only
when $y_{i-1}=y_{i+1}=0$, the number of $2$-cycles in the disjoint
cycle representation of $\Wolf_i^{(54)}$ is once again $2^{n-3}$,
which is even for $n>4$, and $G_{54}$ lies in $\alt_{2^n-1}$.

\begin{conj}[Rule $54$]\label{conj:54}
  For all $n>4$, $G_{54} \cong \alt_{2^n-1}$.
\end{conj}

We have verified Conjecture~\ref{conj:54} explicitly up to $n=8$ as
well.  The reason Conjecture~\ref{conj:54} is difficult to establish
abstractly is that, as with rule $57$, the permutations
$\Wolf_i^{(54)}$ are far removed from the standard generating sets of
$\alt_{2^n-1}$.

%%%%%%%%%%%%%%%%%%%%%%%%%%%%%%%%%%%%%%%%%%%%%%%%%%%%%%
\section{Seven exceptions}\label{sec:exceptional}
%%%%%%%%%%%%%%%%%%%%%%%%%%%%%%%%%%%%%%%%%%%%%%%%%%%%%%

Of the $41$ representative Wolfram rules, we have seen that $26$ have
trivial dynamics groups, $9$ are invertible, and $1$ rule is both.  In
this section we study the remaining $7$ rules.  One key to
understanding the dynamics of these noninvertible rules with
nontrivial dynamics is to note that they agree with an invertible rule
on a union of its orbits, which implies that the dynamics group of the
noninvertible rule is a homomorphic image of the dynamics group of the
invertible rule with which it agrees.  We begin by establishing these
general facts.

\begin{thm}[Independence and invertibility]\label{thm:ind-inv}
  If $\Fy$ is a $\pi$-independent \sds, then there exists an
  invertible \sds $\Fy'$ such that (1) $\Fy$ and $\Fy'$ agree on the
  set $\per(\Fy)$, (2) $\per(\Fy)$ is a union of orbits of the action
  of $\dg(\Fy')$ on $\F^n$, and (3) there is a surjection from
  $\dg(\Fy')$ onto $\dg(\Fy)$.
\end{thm}

\begin{proof}
  We define new local functions $F_i'$ as follows.  Call a state
  $\y\in \F^n$ \emph{$i$-periodic} if the sequence $\y$, $F_i(\y)$,
  $F_i(F_i(\y))$, $F_i(F_i(F_i(\y)))$, etc., eventually returns to
  $\y$; otherwise call it \emph{$i$-transitory}.  The function $F_i'$
  is defined as equal to $F_i$ on the $i$-periodic states and the
  identity function of the $i$-transitory ones.  It is now
  straightforward to check that $F_i'$ is both $Y$-local at $i$ and
  bijective.  By Definition~\ref{def:invertible}, the bijectivity of
  each $F_i'$ means that $\Fy'$ is invertible.  Next, by
  Proposition~\ref{prop:bij} each $F_i^*$ is a permutation of
  $\per(\Fy)$. Thus the states in $\per(\Fy)$ are $i$-periodic for
  each $i$, and $\Fy$ and $\Fy'$ agree completely on this set of
  states.  The invariance of $\per(\Fy)$ under each $F_i$ and its
  agreement with the corresponding $F_i'$ means that $\per(\Fy)$ is a
  union of orbits under the action of $\dg(\Fy')$.  And finally, the
  projection from the subdirect product containing $\dg(\Fy')$ onto the
  factors corresponding to the orbits contained in $\per(\Fy)$
  produces the required surjective group homomorphism.
\end{proof}

\begin{cor}[Independence and invertibility]\label{cor:ind-inv}
  If $\Fy$ is a $\pi$-independent \sds with only two possible
  vertex states, then the update rules for an invertible \sds $\Fy'$
  satisfying the conclusions of Theorem~\ref{thm:ind-inv} can be
  obtained by replacing every {\tt 0} and {\tt 1} in the tags of the
  update rules of $\Fy$ with the symbol {\tt -}.
\end{cor}

\begin{proof}
  The assertion merely describes the effect the previous construction
  has on tags.
\end{proof}

As an illustration, we apply the corollary to each of the $7$
exceptional rules.

\begin{exmp}[Exceptional rules and invertibility]
  By Corollary~\ref{cor:ind-inv}, rules $28$~({\tt 0x--}) and
  $29$~({\tt 0x-1}) agree with rule $156$ ({\tt -x--}) when restricted
  to $P_{28}$ and $P_{29}$, respectively.  Similarly, rules $1$~({\tt
    000x}), $9$~({\tt 00-x}), $129$~({\tt -00x}), $137$~({\tt -0-x})
  and $73$~({\tt 0--x}) agree with rule $201$~({\tt ---x}) when
  restricted to $P_{1}$, $P_{9}$, $P_{129}$, $P_{137}$, and $P_{73}$,
  respectively.
\end{exmp}

The proof of Theorem~\ref{thm:ind-inv} also provides a concrete
description of the periodic states.

\begin{cor}[Periodic orbits]\label{cor:orb-ind-inv}
  If $\Fy$ is a $\pi$-independent \sds with only two possible
  vertex states and $\Fy'$ is the corresponding invertible \sds
  satisfying the conclusions of Theorem~\ref{thm:ind-inv}, then the
  nonperiodic states under $\Fy$ are the union of the orbits under
  $\Fy'$ that contain a state $\y$ that can be altered as a
  consequence of a {\tt 0} or a {\tt 1} in the tag of an update rule
  of $\Fy$.
\end{cor}

Table~\ref{tbl:exceptional} contains a summary of our results about
these $7$ exceptional rules.  These noninvertible rules with
nontrivial dynamics naturally fall into three classes: $\{28,29\}$,
$\{1,9,129,137\}$, and $\{73\}$.

\begin{table}
\begin{tabular}{|c|c|rrr|c|c|c|c|}
  \hline
 Rule & tag & inv & refl& i.r.& $o_k$ & $P_k$ & $c_k$ & $G_k$ \\ \hline\hline
   28 & \tt{0x--} & 199 &  70 & 157 & $\card{\no_{AB}}+1$ & $\no_{BCF}\cup\zer$ & 2 & $\Z_2^n$ \\
   29 & \tt{0x-1} &  71 &  71 &  29 & $\card{\no_{AB}}$ & $\no_{BCF}$ & 2 & $\Z_2^n$ \\ \hline
    1 & \tt{000x} & 127 &   1 & 127 & $1$ & $\no_A$ & 6 & $\lucgp_n$ \\
    9 & \tt{00-x} & 111 &  65 & 125 & $1$ & $\no_A$ & 6 & $\lucgp_n$ \\
  129 & \tt{-00x} & 126 & 129 & 126 & $2$ & $\no_A\cup\one$ & 6 & $\lucgp_n$  \\
  137 & \tt{-0-x} & 110 & 193 & 124 & $2$ & $\no_A\cup\one$ & 6 & $\lucgp_n$  \\ \hline
   73 & \tt{0--x} & 109 &  73 & 109 & $\card{\no_{CD}} $ & $\no_C$ & 6 & Thm~\ref{thm:73} \\ \hline
\end{tabular} \vspace*{2mm}
\caption{The seven exceptional cases.\label{tbl:exceptional}}
\end{table}

%%%%%%%%%%%%%%%%%%%%%%%%%%%%%%%%
\subsection*{Rules $28$ and $29$}

By Corollary~\ref{cor:ind-inv}, both $G_{28}$ and $G_{29}$ are
quotients of $G_{156}$.  The first step is to decide which states are
in $P_{28}$ and $P_{29}$.  By Corollary~\ref{cor:orb-ind-inv}, rule
$28$ (with tag {\tt 0x--}) must exclude all orbits that contain a
state with subword $111$ and rule $29$ (with tag {\tt 0x-1}) must
exclude all orbits that contain a state with subword $111$ or $000$.
For rule $29$, we are left with those states that have no $10$ blocks
of length $4$ or more (because any orbit that contains $1000$ or
$1100$ also contains $1110$) and for rule $28$ we are left with the
same states plus $\zer$ which is fixed.  These states can be described
by excluding $B$ = `111', $C$ =`000' and $F$=`1100'.  Because fixed
states contribute trivial groups to the subdirect product and both
$G_{28}$ and $G_{29}$ are obtained by restricting the action of
$G_{156}$ on $\F_2^n$ to the same set of non-fixed periodic states,
$G_{28}$ and $G_{29}$ are isomorphic groups.  The number of orbits
depends on the number of ways to tile $n$ cyclically arranged
positions by tiles of size $2$ and $3$ and, as we saw in
Section~\ref{sec:trivial}, this number is counted by $\card{N_{AB}}$.
Thus $o_{29} = \card{N_{AB}}$ and $o_{28} = \card{N_{AB}}+1$.
Continuing the notation used in our discussion of invertible rule
$156$, this means that $G_{28}=G_{29}$ is the image of $G_{156}$ in
the product $\sym_2^{(1)} \times \sym_2^{(2)} \times \cdots \times
\sym_2^{(n)}$ which is isomorphic to $\Z_2^n$.  And since $G_{156}$
projects fully onto each factor, $G_{28}=G_{29}$ is all of $\Z_2^n$.
As a final note, when $n=4$, $m+1$ is not allowed to be $3$ and the
factors $\sym_2^{\ell}$ do not occur in the subdirect product.  Thus
$G_{28} =G_{29}$ are trivial groups for $n=4$.

An alternative way to proceed would have been to calculate
$c_{28}=c_{29}=2$, conclude that $G_{28}=G_{29}$ is a quotient of
$\Z_2^n$ and then to argue that it is actually equal to $\Z_2^n$.
Finally, we remark that rules $28$ and $29$ have the same dynamics
group as rule $51$ but for very different reasons.  It is perhaps
surprising that rules $28$ and $29$ have such simple dynamics groups
since they were among the six exceptional cases that needed to be
dealt with separately in the proof of Theorem~\ref{thm:104} due to the
complications in analyzing their dynamics.

%%%%%%%%%%%%%%%%%%%%%%%%%%%%%%%%%%%%%%%%%%%%
\subsection*{Rules $1$, $9$, $129$ and $137$}

By Corollary~\ref{cor:ind-inv}, $G_1$, $G_9$, $G_{129}$ and $G_{137}$
are all quotients of $G_{201}$.  To decide which quotient, the first
step is to decide which states are in $P_1$, $P_9$, $P_{129}$ and
$P_{137}$.  By Corollary~\ref{cor:orb-ind-inv}, all four sets must
exclude any state with adjacent ones (although $P_{129}$ and $P_{137}$
include the fixed state $\one$) and, by our previous analysis of rule
$201$, the states with only isolated $1$s form a single orbit under its
action and thus are included in all four periodic sets.  This means
that $o_1 = o_9 = 1$ and $o_{129} = o_{137} = 2$.  Moreover, recalling
that fixed states contribute trivial groups to the subdirect product,
$G_1$, $G_9$, $G_{129}$ and $G_{137}$ are all isomorphic to the group
obtained when the action of $G_{201}$ is restricted to the single
orbit $\no_A$.  In other words, the groups $G_1$, $G_9$, $G_{129}$,
and $G_{137}$ are all isomorphic to $\lucgp_n$.

%%%%%%%%%%%%%%%%%%%%%%%
\subsection*{Rule $73$}

Our final rule, rule $73$, is also a quotient of $G_{201}$ as
indicated by Corollary~\ref{cor:ind-inv}, but this time the number of
orbits retained is much larger.  In fact, by
Corollary~\ref{cor:orb-ind-inv}, only states containing $C$=`111' as a
subword need to be removed.  Since neither rule $201$ nor rule $73$
can bring subwords of the form `111' into existence, all states in
$\no_C$ are periodic under $G_{73}$.  Because the orbits under rule
$201$ each contain a unique state with no isolated $1$s, the number of
orbits in $P_{73}$ equals the number of states with no `111' and no
`010'.  Thus $P_{73}=\no_C$ and $o_{73}=\card{\no_{CD}}$.  When the
orbit groups for these orbits are factored into direct products and
redundancies are removed, the Lucas group $\lucgp_n$ and all of the
Fibonacci groups $\fibgp_m^{(\ell)}$ arise, except those with $m=n-3$
and $m=n-4$.  This is because the complement of the $0$-block of
length $m$ starting at location $\ell$ must be tilable with pairs of
adjacent $1$s that alternate with other $0$-blocks.  When $m=n-2$ the
complement is a single pair of $1$s and when $m \leq n-5$ the
complement can be a word of the form $1^20^i1^2$ with $i>0$.

\begin{thm}\label{thm:73}
  The group $G_{73}$ is the projection of $G_{201}$ into the following
  direct product:
  \[ G_{73} \subset \lucgp_n \times \prod_{\ell=1}^{n} \left(
    \fibgp_{n-2}^{(\ell)} \times \prod_{m=3}^{n-5}
    \fibgp_m^{(\ell)}\right) \]
\end{thm}

As with rule $201$, we used a computer program to calculate the order
of $G_{73}$ for $4\leq n\leq 7$, and the results are shown in
Table~\ref{tbl:73-201-156}.

%%%%%%%%%%%%%%%%%%%%%%%%%%%%%%%%%%%%%%%%%%%%%%%%%%
\section{Concluding Remarks}\label{sec:conclusion}
%%%%%%%%%%%%%%%%%%%%%%%%%%%%%%%%%%%%%%%%%%%%%%%%%%

Even though the periodic states and the dynamics groups have now been
largely described for all $\pi$-independent \acas, many interesting
research topics remain. Examples include proving or disproving
Conjectures~\ref{conj:fg-lg}, \ref{conj:57}, and~\ref{conj:54}.
The theory and techniques developed in this paper also apply to the
larger class of $\pi$-independent \sdss. In contrast to \acas, where
the graph is $\Circ_n$, \sdss are defined over arbitrary finite
graphs.  Our initial work in~\cite{Macauley:08a} started from \acas
since this class of systems is amenable to analysis and still exhibits
interesting behavior. Extending our analysis to \sdss over arbitrary
graphs poses a bigger challenge, even when restricted to special
classes of vertex functions such as logical NOR and NAND functions
(which are always $\pi$-independent~\cite{Hansson:05b}), or those
inducing invertible \sdss.
As may be clear from this paper, determining the dynamics group from
its definition may be challenging. It would be interesting to
investigate whether there is a result similar to the celebrated
Seifert--van Kampen Theorem that would allow one to deduce dynamics
groups based on, e.g., graph unions or minors. Even if this may be too
much to expect in the general case, it would still be interesting even
if it applies for special classes of vertex functions.
Another problem that may be worth pursuing in the general context is
how to give an efficient presentation of the dynamics groups. In
particular, they are all finite quotients of Coxeter groups. What are
the additional relations arising from the fact that vertex functions
are defined using $\Z_2$ as the state space?

The original motivation for this work was to explore the concept of
$\pi$-independent \acas and the possible dynamics groups that can
arise. Admittedly, this theory and associated techniques still need to
be developed in order for this to become a powerful tool in the study
of discrete dynamical systems. Nonetheless, the construction of the
dynamics group establishes a new connection between algebra and
discrete dynamical systems. As such, it provides a possible avenue for
extending \sds theory through a large body of established mathematical
theory.

%% ------------------------------------------------------------
\medskip\noindent
\subsection*{Acknowledgments}

The second author gratefully acknowledges the support of the National
Science Foundation. The first and third authors thank the Network
Dynamics and Simulation Science Laboratory (NDSSL) at Virginia Tech
for the support of this research which has been partially supported by
HSD Grant SES-0729441, CDC Center of Excellence in Public Health
Informatics Grant 2506055-01, NIH-NIGMS MIDAS project 5 U01
GM070694-05, and DTRA CNIMS Grant HDTRA1-07-C-0113.

%% ------------------------------------------------------------

\bibliographystyle{amsplain}

%%%%%%%%%%%%%%
\end{document}